# Simultaneous model identification and optimization in presence of model-plant mismatch

Jasdeep S. Mandur, Hector M. Budman*

Department of Chemical Engineering, University of Waterloo, Waterloo, ON, N2L3G1, Canada

**Abstract**

In a standard optimization approach, the underlying process model is first identified at a given set of operating conditions and this updated model is, then, used to calculate the optimal conditions for the process. This "two-step" procedure can be repeated iteratively by conducting new experiments at optimal operating conditions, based on previous iterations, followed by re-identification and re-optimization until convergence is reached. However, when there is a model-plant mismatch, the set of parameter estimates that minimizes the prediction error in the identification problem may not predict the gradients of the optimization objective accurately. As a result, convergence of the "two-step" iterative approach to a process optimum cannot be guaranteed. This paper presents a new methodology where the model outputs are corrected explicitly for the mismatch such that, with the updated parameter estimates the identification and optimization objectives are properly reconciled. With the proposed corrections being progressively integrated over the iterations, the algorithm has guaranteed convergence to the process optimum and also, upon convergence, the final corrected model predicts the process behavior accurately. The proposed methodology is illustrated in a run-to-run optimization framework with a fed-batch bioprocess as a case study.

**Highlights**

- An iterative optimization algorithm is proposed.
- The model outputs are corrected iteratively to account for model-plant mismatch.
- Parameters are estimated to satisfy both identification and optimization objectives.
- The proposed algorithm is illustrated using a fed-batch bioprocess.

*Corresponding Author. Tel: +1-519-888-4567 x36980; Fax: +1-519-746-4979; Email: hbudman@uwaterloo.ca

## 1. Introduction

Mathematical models have become an integral part of the process development and subsequent production environment. Besides providing novel insights into the underlying process, they are also used in various model-based optimization and optimal control strategies. When a model is an exact representation of the actual process and is calibrated against noise-free process data, optimizing the model is identical to optimizing the process itself. In such case, the optimal policies derived from model-based optimization can be applied in an "open loop" fashion. However, the above conditions are extremely difficult to meet in practical situations. In the presence of any model uncertainty resulting from either incorrect model structure or measurement noise, the model-based optimization algorithms will result in sub-optimal policies or, in a worst case scenario, may also result in violation of process constraints. To tackle this problem, one possible solution is to search for optimal policies that are robust to model uncertainties (Beyer et al., 2007; Samsatli et al., 1998; Diwekar et al, 1996; Nagy et al., 2004; Ruppen et al., 1995; Terwiesch et al., 1994). Although this approach can ensure feasibility within *a priori* known bounds of uncertainty, the optimal policies are often conservative and, in some processes, may lead to significant loss in economic objectives. As an alternative, another possibility is to use an iterative approach where the model is updated using new measurements at previously calculated optimal policy and the updated model is, then, re-optimized for the next optimal policy (Ruppen et al., 1998; Eaton et al., 1990; Chen et al., 1987; Marlin et al., 1997). This process is referred to as a "two-step" approach and is repeated until a convergence is achieved.

This paper deals with the application of latter approach to batch/fed-batch processes. Assuming the process data is available only at the end of batch, the problem is solved in a run-to-run optimization framework. However, there is no restriction on applying the proposed algorithm to online optimization problems if measurements are available online.

The convergence of the standard two-step approach is governed by; (1) whether the sub-optimal policies provide enough excitation to update all the parameters and (2) how close the model can describe the actual process. The first condition can be addressed, to a certain extent, by incorporating design of experiments in the optimization objectives (Martinez et al., 2013). In this way, a trade-off between the optimal policies and the policies that generate more informative process data can be achieved. Then, if the model is a true of representation of the process, the two-step approach will converge to the actual process optimum, where the total number of iterations needed for convergence, will depend on measurement noise and the extent of excitation.

Regarding the second condition, mentioned above, model-plant mismatch is inevitable in almost all practical applications. In an attempt to capture the process behavior accurately, models often become too complex and computationally demanding rendering them unsuitable for optimization or control purposes. Also, with a limited number of measurable states and in the presence of measurement noise, it is not possible to estimate all the parameters accurately and, therefore, model reduction techniques are often used to reduce the number of parameters that can be identified from the given process data. Because of these reasons, one has to generally rely on simpler but inaccurate model structures for optimization and control. If the structural inaccuracy is not considered explicitly in the model, calibrating the model over different operating conditions may result in significantly different parameter estimates in order to compensate for the model error around different operating points. As a result of this parametric variability, it is possible that the optimization objectives may get compromised. The change in parameter estimates may be of such an extent that the predicted gradients of the optimization objectives no longer coincide with the gradients measured from the process or, in a worst case scenario, they may even get reversed thus leading to premature convergence to sub-optimal operating policies.

For convergence to the process optimum, it is necessary that the model accurately predicts the optimality conditions of the process as given by the first-order Karush-Kuhn-Tucker (KKT) conditions. Following this idea, a class of algorithms has been developed where the optimization objectives and the constraints are corrected for the bias as well as the difference between their predicted and measured gradients (Roberts et al., 1979; Tatjewski, 2002; Gao et al., 2005; Chachuat et al., 2009; Marchetti et al., 2009; Costello et al., 2011). These algorithms differ in the way the model is updated and on how the modifications to the objective function and constraints are implemented. For instance, in their pioneer work, Roberts et al. (1979) modified only the optimization objective to account for the difference between predicted and measured output derivatives assuming the constraints are process independent and are known. The modification term was based on the Lagrange multipliers where the Lagrangian function was obtained by integrating the identification and optimization objectives. In subsequent studies, Tatjewski (2002) and Gao et al. (2005) replaced the parameter estimation step by introducing a linear term in the outputs that corrects for the difference between the predicted and measured outputs. In a more recent version of these algorithms, referred to as modifier adaptation (Chachuat et al., 2009; Marchetti et al., 2009 and Costello et al., 2011), the authors eliminated the model update step altogether and updated only the optimization objectives based on differences between the gradients in the optimality condition. Since the convergence to a process optimum was driven solely by the correction in the optimization gradients, the final model-based optimal solution remained unaffected by the elimination of the model-update step. However, this approach results in a bias between the predicted and measured outputs and as result, the algorithm can no longer be applied to the problems where prediction accuracy is required. One such case, as recently pointed out by Costello et al. (2011), is when the optimal input profiles are implemented within a closed-loop control to ensure that the process is operated to meet safety and environmental constraints. Here, the model is required to provide accurate reference trajectories for low-level controllers. The prediction accuracy of the model is also very relevant for biotechnological processes where it is important to predict the evolution of toxic by-products along a batch culture. Thus, to address a broad range of problems, it is very important to satisfy both identification and optimization objectives at the optimum. One of the major bottlenecks in implementation of this class of algorithms is their sensitivity to the noise in measured gradients (Marchetti et al., 2009). To avoid too much aggressiveness in the corrections and to achieve a smoother convergence, the corrections have to be filtered using an empirical filter.

In another class of algorithm, Srinivasan et al. (2002) proposed an alternate approach where the identification objective is modified to account for the difference between predicted and measured optimality conditions. By this modification, the parameter estimates can be obtained so as to provide a trade-off between the identification and optimization objectives based on preselected weights.

In this work, we proposed a linear correction to the model outputs in a way that the updated model parameters not only minimize the bias between the predicted and measured outputs, but it also correct for the optimization objectives. The corrections made over the previous iterations are progressively integrated and by implementing this progressive correction in the model, the conflict between the identification and optimization objectives is reduced significantly. Another key advantage of this approach is that it provides a model-based filter that is shown to outperform the external exponential filter, used in the previous studies, in terms of the rate at which convergence can be achieved. A preliminary discussion of this methodology has been presented in Mandur et al. (2013).

The contents of this paper are organized as follows: Section 2 provides a brief background on two-step approach and modifier adaptation algorithms and also discusses the motivation in detail. Section 3 presents the methodology and theory behind the proposed model correction. The methodology is then illustrated with a case study in Section 4 and finally Section 5 concludes the paper.

## 2. Preliminaries

Let us consider a process model, described by a set of differential equations as follows:

$$\dot{\mathbf{x}} = f(\mathbf{x}, \boldsymbol{\theta}, \mathbf{u}, \mathbf{t}) + \boldsymbol{\upsilon}$$

$$\mathbf{y}_m = h(\mathbf{x}) + \boldsymbol{\eta} \qquad (P.1)$$

Where, $\mathbf{x} \in \mathbb{R}^{n_x}$ is the vector of model states, $\boldsymbol{\theta} \in \mathbb{R}^{n_\theta}$ is the vector of model parameters, $\mathbf{u} \in \mathbb{R}^{n_u}$ is the vector of inputs, $\mathbf{y}_m \in \mathbb{R}^{n_y}$ is the vector of measured output variables, $f \in \mathbb{R}^{n_x}$ is a set of differential equations based on mass and energy balances, $h \in \mathbb{R}^{n_y}$ is a mapping between the model states and predicted outputs and $\boldsymbol{\upsilon}$ and $\boldsymbol{\eta}$ are the vectors of uncertainties representing modelling and measurement errors respectively.

The standard two-step optimization approach starts with a model identification step where the model is calibrated using the process measurements at some initial input conditions. The identification objective is generally based on the minimization of the errors between predicted and measured outputs. For example, the standard least squares estimation problem can be formulated as follows:

$$\theta_k = \arg\min_{\theta} \sum_{i=1}^{N} \|y_m(u_k) - y(u_k, \theta)\|^2$$

$$\text{s.t.} \quad \dot{x} = f(x, \theta, u_k, t)$$

$$y = h(x) \tag{P.2}$$

Where, $y \in \mathbb{R}^{n_y}$ is the vector of predicted outputs, N is the number of time points and subscript $k$ is the current iteration. The identification is then followed by an optimization step, formulated as follows:

$$u_{k+1} = \arg\min_{u} \phi(y, u, \theta_k)$$

$$\text{s.t.} \quad \dot{x} = f(x, \theta_k, u, t)$$

$$y = h(x)$$

$$g(y, u, \theta_k) \leq 0 \tag{P.3}$$

Where, $\phi$ is the objective function to be minimized and $g \in \mathbb{R}^{n_g}$ is a vector of equalities or inequalities.

Let the functions $\phi$ and $g_j$ be continuously differentiable at a set of input conditions $u^*$. If $u^*$ is a process optimum, then there exists a unique vector, $\mu \in \mathbb{R}^{n_g}$ such that:

$$\frac{\partial \phi(y_m, u^*)}{\partial u} + \mu^T \frac{\partial g(y_m, u^*)}{\partial u} = 0 \tag{1a}$$

$$\mu^T g(y_m, u^*) = 0 \tag{1b}$$

$$\mu \geq 0 \tag{1c}$$

$$g(y_m, u^*) \leq 0 \tag{1d}$$

These conditions are collectively known as Karush-Kuhn-Tucker (KKT) conditions, where $\mu$ is a vector of KKT multipliers. For the model-based optimal solution to converge to $u^*$, it is necessary that the model predicts the KKT above conditions accurately. Since $u^*$ is not known a priori, this can be only guaranteed if the identification step (P.2) results in a unique set of model parameters ($\theta_k$) such that the model satisfies the following conditions for all set of values of $u \in \mathbb{R}^{n_u}$:

$$\frac{\partial \phi(y, u, \theta_k)}{\partial u_i} = \frac{\partial \phi(y_m, u)}{\partial u_i} \tag{2a}$$

$$\frac{\partial g_j(y, u, \theta_k)}{\partial u_i} = \frac{\partial g_j(y_m, u)}{\partial u_i} \tag{2b}$$

When there is only measurement noise (i.e. $\upsilon = 0$), the above conditions can be satisfied over a finite number of iterations where:

$$\frac{\sum_{k=1}^{N} \theta_k}{N} \sim \theta^* \tag{3}$$

Where, $\theta^* \in \mathbb{R}^{n_\theta}$ is the set of parameter values satisfying the conditions (2a) and (2b)

However, in the presence of model structure error ($\upsilon \neq 0$), $\theta^*$ does not exist. Since the error term $\upsilon$ represents the unmodelled dynamics of the process, it is a time varying function of model states ($x$) and inputs ($u$). Then, when the inaccurate model, given by $\dot{x} = f(x, \theta, u_k, t)$, is calibrated over different input conditions, the parameter estimates change so as to compensate for the modelling error which varies with respect to the input conditions. Consequently, there is no unique set of parameter estimates that can satisfy the identification objective (P.2) over the entire space of input conditions. For any particular set of parameter estimates ($\boldsymbol{\theta}_k$), the model is accurate only in the neighbourhood of the corresponding input values ($\mathbf{u}_k$) whereas away from this region, the prediction accuracy of the model continues to decrease as the distance increases. As a result, the model may not predict the gradients of the optimization objective and constraints accurately. In a worst case scenario, it is also possible that the change in model parameters as the input conditions change is of such an extent that the predicted and measured gradients have opposite signs, in which case the model-based optimization can no longer drive the changes in the inputs in the direction of process optimum. This implies:

$$\phi(\mathbf{y_m}, \mathbf{u_{k+1}}) > \phi(\mathbf{y_m}, \mathbf{u_k}) \tag{4}$$

In this case, the two-step approach will converge to a non-optimal solution where the measured gradients are still non-zero, or in other words, the predicted KKT conditions do not match with those measured from the process. Therefore, to ensure convergence of the algorithm to a process optimum, the differences between the predicted and measured gradients of the optimization problem must be eliminated or at least minimized at each intermediate input condition as shown schematically in Figure 1.

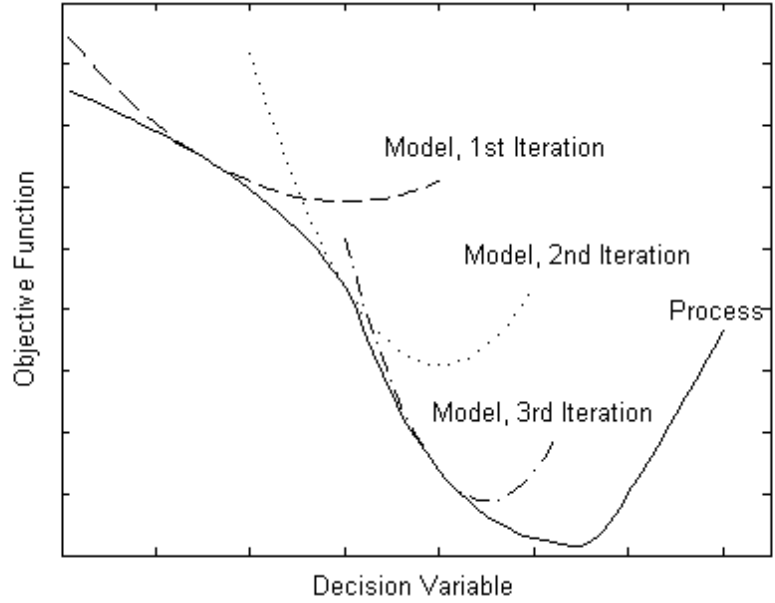

Figure 1: Iterative improvements in model-based optimal solution

As stated in the Introduction, *modifier adaptation* algorithms enforce the matching conditions (Equations (2a) and (2b)) by adding correction terms directly to the corresponding optimization quantities. Accordingly, the modified optimization problem is as follows:

$$\mathbf{u_{k+1}} = \arg \min_{\mathbf{u}} \phi(\mathbf{y}, \mathbf{u}, \boldsymbol{\theta}) + \boldsymbol{\lambda}_{\phi k}^{\mathbf{T}} \mathbf{u}$$

$$\text{s.t.} \quad \dot{\mathbf{x}} = f(\mathbf{x}, \boldsymbol{\theta}, \mathbf{u}, \mathbf{t})$$

$$\mathbf{y} = h(\mathbf{x})$$

$$\mathbf{g}(\mathbf{y}, \mathbf{u}, \boldsymbol{\theta}) + \boldsymbol{\varepsilon}_{gk} + \boldsymbol{\lambda}_{gk}^{\mathbf{T}}(\mathbf{u} - \mathbf{u_k}) \leq \mathbf{0} \quad \quad (\text{P.4})$$

Where, $\boldsymbol{\lambda}_\phi$ and $\boldsymbol{\lambda}_g$ are referred to as modifiers that are used to correct for the gradients of objective function and constraints respectively and $\boldsymbol{\varepsilon}_g$ is a modifier introduced to correct for the bias in predicted and measured constraints. The corrections are calculated at the $k^{th}$ iteration as follows:

$$\lambda_{\phi k_i} = \frac{\partial \phi(\mathbf{y_m}, \mathbf{u_k})}{\partial u_i} - \frac{\partial \phi(\mathbf{y}, \mathbf{u_k}, \boldsymbol{\theta})}{\partial u_i} \qquad (5a)$$

$$\lambda_{g k_{ij}} = \frac{\partial g_j(\mathbf{y_m}, \mathbf{u_k})}{\partial u_i} - \frac{\partial g_j(\mathbf{y}, \mathbf{u_k}, \boldsymbol{\theta})}{\partial u_i} \qquad (5b)$$

$$\varepsilon_{g_j} = g_j(\mathbf{y_m}, \mathbf{u_k}) - g_j(\mathbf{y}, \mathbf{u_k}, \boldsymbol{\theta}) \qquad (5c)$$

To avoid excessive corrective actions and to reduce the sensitivity to measurement noise, these corrections are filtered before implemented in (P.4) as follows:

$$\boldsymbol{\Lambda}_k = \mathbf{K}\boldsymbol{\Lambda}'_k + (1 - \mathbf{K})\boldsymbol{\Lambda}_{k-1} \qquad (6)$$

Where, $\boldsymbol{\Lambda}$ represent the vector of modifiers defined as: $\boldsymbol{\Lambda} = [\boldsymbol{\lambda}_\phi, \boldsymbol{\lambda}_g, \boldsymbol{\varepsilon}_g]$ and $\mathbf{K}$ represents the filter gain.

It should be noted that whether the model is updated or not, the corrected optimization objective does not depend on the model. Therefore, the aggressive changes in the inputs, resulting from inaccurate predictions of the modified objective function or constraints, have to be controlled by the filter gain which may further reduce the speed of convergence. Also, there is no systematic way to choose a priori the filter gain. In the case study presented later, it is observed that when the model is used to correct for the errors in gradients as proposed in the current study, it results in more accurate predictions during the search for optimal solution thus leading to faster convergence with less oscillatory behaviour.

## 3. Proposed Methodology

The basic idea in the proposed methodology is to search for model parameters such that the differences between the predicted and measured gradients of optimization problem given in Equations (2a) and (2b) are minimized along with the minimum prediction error from identification problem (P.2). However, since the identification and optimization objectives are independent of each other, in the presence of model structure error, the values of model parameters that satisfy both objectives do not exist. In other words, the parameter estimates that minimize the difference in gradients may not minimize the prediction error at the same time. To this end, a linear correction term is added to the model outputs such that the conflicting objectives can be reconciled.

Let $\epsilon_k$ be the minimum sum of squared errors between the predicted and measured outputs in problem (P.2), corresponding to the parameter estimates $\boldsymbol{\theta_k}$ as follows:

$$\epsilon_k = \sum_{i=1}^{N} \|\mathbf{y_m}(\mathbf{u_k}) - \mathbf{y}(\mathbf{u_k}, \boldsymbol{\theta_k})\|^2 \tag{7}$$

Let $\Delta\boldsymbol{\theta_k}$ be the change in parameter estimates, with respect to $\boldsymbol{\theta_k}$, required to minimize the difference between the predicted and measured gradients (Equations (2a) and (2b)) at $k^{th}$ iteration. The updated model parameters i.e. $\boldsymbol{\theta_k} + \Delta\boldsymbol{\theta_k}$, then, no longer minimizes the updated sum of squared errors. This is schematically illustrated in Figure 2 which shows a probability density function of a model parameter centered on the estimate $\theta_k$, calculated by least squares.

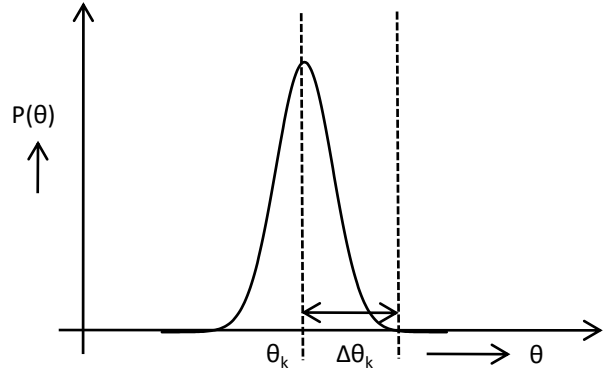

Figure 2: Illustration of lack of fit

The corresponding sum of squared errors using the perturbed parameter value $\boldsymbol{\theta_k} + \Delta\boldsymbol{\theta_k}$ can be represented as:

$$\epsilon'_k = \sum_{i=1}^{N} \|\mathbf{y_m}(\mathbf{u_k}) - \mathbf{y}(\mathbf{u_k}, \boldsymbol{\theta_k} + \Delta\boldsymbol{\theta_k})\|^2 \qquad (8)$$

Let us introduce a vector of corrections $-\mathbf{c_k}$ to the model outputs such that the sum of squared errors, $\epsilon'_k$ with the corrected model remains equal to the original value of $\epsilon_k$ (given by Equation (7)):

$$\sum_{i=1}^{N} \|\mathbf{y_m}(\mathbf{u_k}) - (\mathbf{y}(\mathbf{u_k}, \boldsymbol{\theta_k} + \Delta\boldsymbol{\theta_k}) - \mathbf{c_k})\|^2 = \sum_{i=1}^{N} \|\mathbf{y_m}(\mathbf{u_k}) - \mathbf{y}(\mathbf{u_k}, \boldsymbol{\theta_k})\|^2 \qquad (9)$$

To satisfy (9), the equality can be satisfied term by term as follows:

$$\mathbf{y_m}(\mathbf{u_k}) - \mathbf{y}(\mathbf{u_k}, \boldsymbol{\theta_k} + \Delta\boldsymbol{\theta_k}) + \mathbf{c_k} = \mathbf{y_m}(\mathbf{u_k}) - \mathbf{y}(\mathbf{u_k}, \boldsymbol{\theta_k}) \qquad (10)$$

The correction term can then be solved from (10) by:

$$c_k = y(u_k, \theta_k + \Delta\theta_k) - y(u_k, \theta_k) \qquad (11)$$

Let the model outputs $y(u_k, \theta_k + \Delta\theta_k)$ be approximated around $\theta_k$ using Taylor Series Expansions. This will result in the following expression:

$$y(u_k, \theta_k + \Delta\theta_k) = y(u_k, \theta_k) + Dy(\theta_k)\Delta\theta_k + \cdots \qquad (12)$$

Where, $D$ is the Jacobian matrix of output derivatives with respect to model parameters. After substituting Equation (12) into (11), the correction term $c_k$ can be expressed as:

$$c_k = Dy(\theta_k)\Delta\theta_k + \cdots \qquad (13)$$

Assuming the model to be linear in the neighborhood of $\theta_k$, the correction term $c_k$ is approximated by the first-order derivative as follows:

$$c_k = Dy(\theta_k)\Delta\theta_k \qquad (14)$$

Since the linear approximation of the model is generally valid only within a certain region around $\theta_k$, for $\Delta\theta_k$ outside this region, it would not be possible to restore the prediction error with the updated model to its minimum as the LHS of the equality, given by Equation (9), may have significant error. Therefore, to enforce the approximate validity of the linear approximation, a constraint on $\Delta\theta_k$ is imposed which is based on a relative truncation error, defined as follows:

$$\epsilon^T = \frac{y(u_k, \theta_k + \Delta\theta_k) - Dy(\theta_k)\Delta\theta_k}{y(u_k, \theta_k)} \qquad (15)$$

Based on its definition, the calculation of $\Delta\boldsymbol{\theta}_k$ is calculated using an optimization problem as follows:

$$\Delta\boldsymbol{\theta}_k = \arg\min_{\Delta\boldsymbol{\theta}} \left( \mathbf{w}_\phi \left| \frac{\partial \phi(\mathbf{y}_m, \mathbf{u}_k)}{\partial \mathbf{u}} - \frac{\partial \phi(\mathbf{y}, \mathbf{u}_k, \boldsymbol{\theta}_k + \Delta\boldsymbol{\theta})}{\partial \mathbf{u}} \right| \right.$$

$$\left. + \mathbf{w}_g \left| \frac{\partial g(\mathbf{y}_m, \mathbf{u}_k)}{\partial \mathbf{u}} - \frac{\partial g(\mathbf{y}, \mathbf{u}_k, \boldsymbol{\theta}_k + \Delta\boldsymbol{\theta})}{\partial \mathbf{u}} \right| \right)$$

s.t. $\quad \dot{\mathbf{x}} = f(\mathbf{x}, \boldsymbol{\theta}_k + \Delta\boldsymbol{\theta}, \mathbf{u}, t)$

$\quad\quad \mathbf{y} = h(\mathbf{x}) - D\mathbf{y}(\boldsymbol{\theta}_k)\,\Delta\boldsymbol{\theta}$

$\quad\quad \boldsymbol{\epsilon}^T \leq \epsilon_{max}^T$ \hfill (P.5)

Where, $\mathbf{w}_\phi$ and $\mathbf{w}_g$ are vectors of normalizing weights for the objective function and constraints gradients respectively and $\epsilon_{max}^T$ is the constraint or limit on truncation error that is imposed to ensure the approximated validity of the linear approximation of the correction term.

It is important to note here that the above constraint on $\Delta\boldsymbol{\theta}_k$ is somewhat equivalent to the filter gain in modifier adaptation algorithms as the restriction on $\Delta\boldsymbol{\theta}_k$ also restricts the ability of the model to predict gradients of the optimization problem exactly which is very critical when the gradients are associated with significant level of noise. On the other hand, in contrast with the filter gain in modifier adaptation algorithms, $\Delta\boldsymbol{\theta}_k$ is based on a physical rationale since it represents an allowable model prediction error.

To summarize the procedure, the estimation of model parameters is divided into two steps:

**Step 1**: The parameters are updated to minimize the error between the outputs as predicted by the previously corrected model and those measured from the process. Let us define this update in parameter estimates by $\Delta \theta_{iden_k}$.

**Step 2:** The change in model parameters $\Delta \theta_k$ and the corresponding model correction is, then, calculated such that the updated model predict the gradients of the optimization problem at current input conditions and at the same time, to adjust the prediction error to the same value obtained in the previous step.

The overall update step can be written as:

$$\theta'_k = \theta'_{k-1} + \Delta \theta_{iden_k} + \Delta \theta_k \qquad (16)$$

The corrected model with the updated parameter estimates $\theta'_k$ is then optimized for the next iteration. It should be noted, here, that the model correction term is being carried forward into the next iteration and as a result, it has a cumulative effect. The prediction inaccuracies of the model continue to decrease as the model is corrected progressively towards the process optimum. Finally, at the optimum, the corrected model simultaneously satisfies both the identification and objective objectives.

3.1. Conditions for Convergence

At a given set of input conditions, let us define a bounded space for model parameters such that $\forall\, \boldsymbol{\theta} \in [\boldsymbol{\theta}_{lb}, \boldsymbol{\theta}_{ub}]$:

$$\text{The corrected model is stable} \tag{17a}$$

$$\nabla^2 \phi(\mathbf{u_k}, \boldsymbol{\theta}'_\mathbf{k}) > 0 \quad \text{(Positive definite)} \tag{17b}$$

Then if the bound on truncation error $\epsilon_{max}^T$ is such that $\forall\, \mathbf{u}_k$:

$$\left| \frac{\partial \phi(\mathbf{y_m}, \mathbf{u_k})}{\partial \mathbf{u}} - \frac{\partial \phi(\mathbf{y}, \mathbf{u_k}, \boldsymbol{\theta}'_\mathbf{k})}{\partial \mathbf{u}} \right| < \varepsilon \tag{C.1}$$

Where, $\varepsilon$ is the tolerance with which the above differences between the measured and predicted gradients of the cost function are minimized, then, the algorithm has a guaranteed convergence towards the process optimum

Since the update in input conditions $\mathbf{u}_k$ is based on model-based optimization, the algorithm will converge only if:

$$\frac{\partial \phi(\mathbf{y}, \mathbf{u_k}, \boldsymbol{\theta}'_\mathbf{k})}{\partial \mathbf{u}} = 0 \tag{18}$$

From condition (C.1), since the predicted gradients are always matches to the ones measured from the process, the Equation (18) holds only when:

$$\left| \frac{\partial \phi(\mathbf{y_m}, \mathbf{u_k})}{\partial \mathbf{u}} \right| < \varepsilon \tag{19}$$

## 3.2. Termination Criteria

Let the algorithm converges to a stationary point $\mathbf{u_k^*}$. Then, at $\mathbf{u_k^*}$:

$$\Delta \boldsymbol{\theta}_{iden_k} = 0 \tag{20a}$$

$$\Delta \boldsymbol{\theta}_k = 0 \tag{20b}$$

Since at convergence $\mathbf{u_k^*} = \mathbf{u_{k-1}}$, the parameter estimates $\boldsymbol{\theta}'_{k-1}$ minimizing the prediction error at $\mathbf{u_{k-1}}$ also minimize the prediction error at $\mathbf{u_k}$ and, therefore, the parameter change in **Step 1** is zero. Similarly no further corrections are required for the gradients as they have already been corrected in the previous iteration. Hence, the update in **Step 2** is also zero.

Equations (20a) and (20b) can then be used to define termination criteria. However, in the presence of measurement noise, the above criteria cannot be exactly achieved. To this end, the convergence of the algorithm can be evaluated in terms of convergence in the probability distribution of model parameters. The difference between the distributions in two successive iterations is calculated using the Kullback-Leibler (K-L) divergence (Cover et al., 1991). If $\mathbf{P_k}$ and $\mathbf{P'_{k-1}}$ are the distributions of parameters $\boldsymbol{\theta}'_{k-1} + \Delta \boldsymbol{\theta}_{k_{iden}}$ and $\boldsymbol{\theta}'_{k-1}$ respectively, the K-L divergence between $\mathbf{P_k}$ and $\mathbf{P'_{k-1}}$ is given by;

$$d(\mathbf{P'_{k-1}} || \mathbf{P_k}) = \int \mathbf{P'_{k-1}}(\boldsymbol{\theta}) \log \frac{\mathbf{P'_{k-1}}(\boldsymbol{\theta})}{\mathbf{P_k}(\boldsymbol{\theta})} d\boldsymbol{\theta} \tag{21}$$

The condition based on Equation (20a) is then formulated as:

$$d(\mathbf{P'_{k-1}} || \mathbf{P_k}) \leq \varepsilon_1 \tag{C.2}$$

Similarly, for Equation (20b), the difference between the distributions of $\boldsymbol{\theta}'_{k-1} + \Delta\boldsymbol{\theta}_{k_{iden}}$ and $\boldsymbol{\theta}'_k$ is measured as $d(\mathbf{P_k} || \mathbf{P'_k})$, where, $P'_k$ is the distribution of $\boldsymbol{\theta}'_k$

$$d(\mathbf{P_k} || \mathbf{P'_k}) \leq \varepsilon_2 \qquad (C.3)$$

3.3. Summary of algorithm

Figure 3 presents the flowchart of the algorithm. The algorithm begins with the identification step using the initial inaccurate model structure. For estimation, the problem posed in (P.2) is used. Then, the update in model parameters and corresponding model correction is calculated to correct for the model for both identification and optimization objectives by using the minimization problem posed in (P.5). The updated model is then optimized for the next operating conditions where the above steps are repeated with the updated model. The procedure is repeated until a termination criteria based on either (20a & 20b) or (C.2 & C.3) are satisfied.

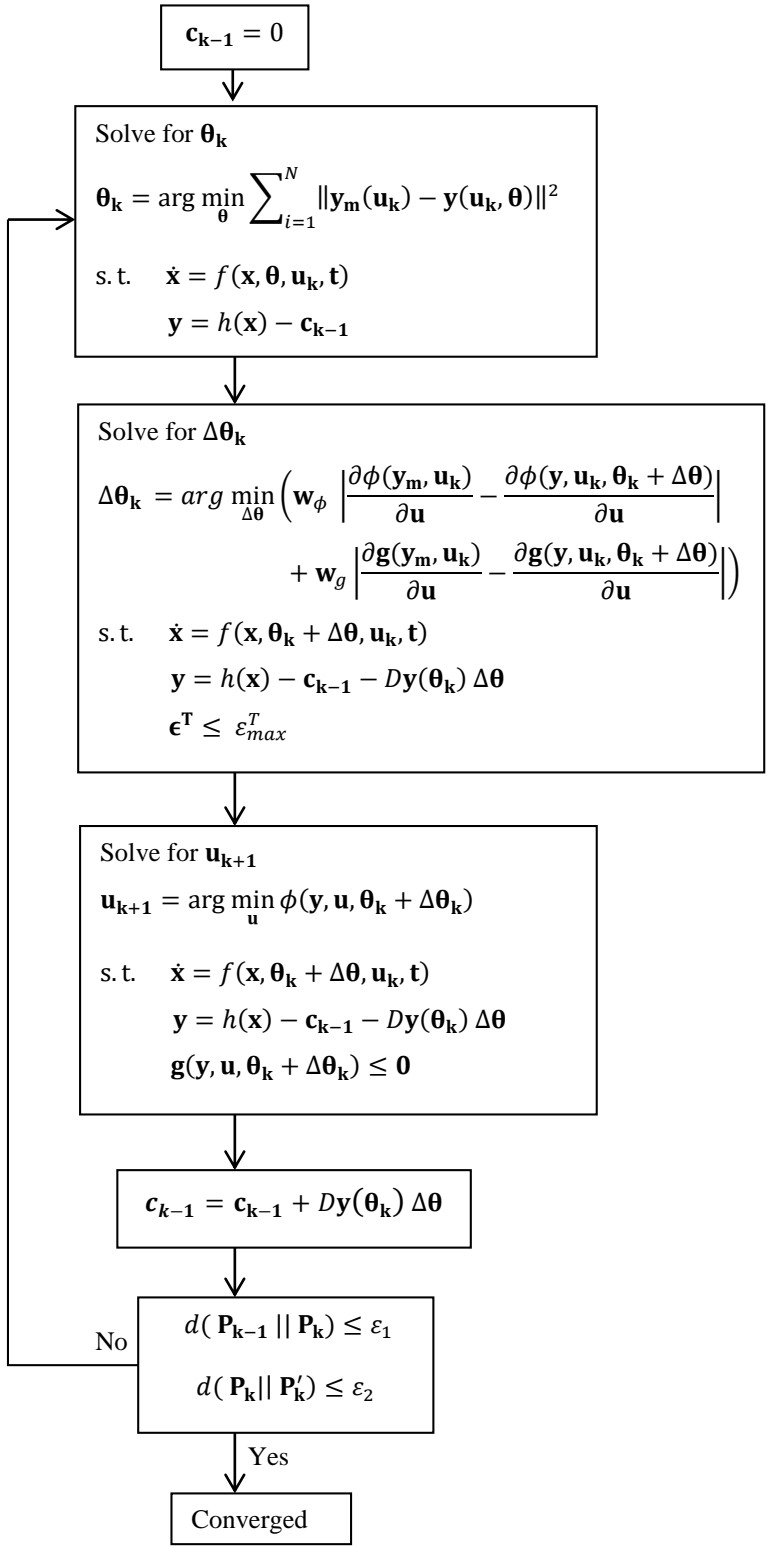

Figure 3. Proposed Algorithm with linear model corrections

## 4. Case Study

4.1. Problem formulation

The proposed optimization algorithm is applied to a penicillin production process where the goal is to maximize the amount of penicillin at the end of batch. To generate the experimental data for model identification and correcting the model for optimization, the *in-silico* experiments are conducted using a process simulator based on the following set of ordinary differential equations (Bajpai and Reuss, 1980; Birol et al. 2002):

$$\frac{dX}{dt} = \left(\frac{\mu_X SX}{K_X X + S}\right) - \frac{X}{V}\frac{dV}{dt} \tag{22}$$

$$\frac{dP}{dt} = \left(\frac{\mu_P SX}{K_P + S + \frac{S^2}{K_I}}\right) - K_H P - \frac{P}{V}\frac{dV}{dt} \tag{23}$$

$$\frac{dS}{dt} = -\frac{1}{Y_{X/S}}\left(\frac{\mu_X SX}{K_X X + S}\right) - \frac{1}{Y_{P/S}}\left(\frac{\mu_P SX}{K_P + S + \frac{S^2}{K_I}}\right) - m_X X + \frac{F s_f}{V} - \frac{S}{V}\frac{dV}{dt} \tag{24}$$

$$\frac{dV}{dt} = F - 6.226 * 10^{-4} V \tag{25}$$

The set of Equations (22)-(24) describes the rate of change in the concentrations of biomass ($X$), penicillin ($P$) and substrate ($S$) respectively and Equation (25) describes the rate of change in the culture volume ($V$). The constants in these equations are defined as follows; $\mu_X$ is the specific growth rate of biomass, $\mu_P$ is the specific rate of penicillin production, $K_X$ and $K_P$ are saturation constants, $K_I$ is a substrate inhibition constant, $K_H$ is a constant representing the rate of consumption of penicillin by hydrolysis, $Y_{X/S}$ and $Y_{P/S}$ are the yields per unit mass of substrate for the biomass and penicillin respectively, $m_X$ represents the consumption rate of substrate needed for maintaining the biomass and $s_f$ is the concentration of substrate in the feed. The values used for these constants are listed in Table 1.

Table 1: Parameters' values for process simulator (Equations 22-25)

| $\mu_x$ | $K_x$ | $\mu_P$ | $K_P$ | $K_I$ | $K_H$ | $Y_{X/S}$ | $Y_{P/S}$ | $m_X$ | $s_f$ |
|---|---|---|---|---|---|---|---|---|---|
| 0.092 | 0.15 | 0.005 | 0.0002 | 0.1 | 0.04 | 0.45 | 0.9 | 0.014 | 600 |

To formulate a model with structural inaccuracy, it is assumed that the user does not have prior knowledge about the consumption of penicillin by hydrolysis and, as a result, the rate of change in the penicillin concentration is modelled as:

$$\frac{dP}{dt} = \left( \frac{\mu_p S X}{K_P + S + \frac{S^2}{K_I}} \right) - \frac{P}{V}\frac{dV}{dt} \qquad (26)$$

Assuming the dynamics for the other states to be known accurately, the uncertain model is then given by the set of Equations (22) and (24-26). To simplify the numerical calculations, it is further assumed that only two model parameters $K_X$ and $K_I$ will be updated in the algorithm whereas the rest of the model parameters are fixed at their nominal values, estimated at initial input conditions as listed in Table 2. The choice of these two parameters as the uncertain ones was based on a preliminary sensitivity analysis.

Table 2: Initial set of input conditions used to estimate the parameters in uncertain model

| | |
|---|---|
| Biomass Conc. ($X_0$) | 0.1 (g/l) |
| Substrate Conc. ($S_0$) | 1 (g/l) |
| Product Conc. ($P_0$) | 0 (g/l) |
| Initial Culture Volume ($V_0$) | 100 (L) |
| Input Feed (F) | 0.04 (L/hr) |

The uncertain model is, then, optimized iteratively as per the procedure summarized in Figure 3, where the final objective is to maximize the amount of penicillin at the end of batch, subject to a terminal constraint on the culture volume. The initial substrate concentration $S_o$ and the input feed rate $F$ are selected as the decision variables whereas the rest of the input variables are fixed at their initial values listed in Table 2. Accordingly, the optimization problem is formulated as follows:

$$\min_{S_o, F} \quad -P(\mathbf{x}, \boldsymbol{\theta}, S_o, F, t_f)$$

$$s.t. \quad (22) \text{ and } (24) - (26)$$

$$V(\mathbf{x}, \boldsymbol{\theta}, S_o, F, t_f) \leq V_{max} \tag{P.6}$$

For reference, the process optimum corresponds to $S_o = 55$ g/l and $F = 0.1728$ L with the final penicillin measured to be $= 592$g

## 4.2. Results and discussion

In the first part of discussion, the convergence properties of the algorithm will be discussed. The performance of the algorithm is evaluated in terms of (1) the rate of convergence and (2) the final converged solution.

The bound on truncation error $\epsilon_{max}^T$ is the major factor that affects the rate of convergence. Let us recall the calculation of $\Delta\boldsymbol{\theta}_k$. The parameter estimates $\boldsymbol{\theta}_k$ that minimize the prediction error at $k^{th}$ iteration may not predict the gradients of the optimization problem correctly for which a change in estimates $\Delta\boldsymbol{\theta}_k$ is required. Then, in order to ensure that the prediction error is also minimized with the updated parameters, the model outputs are corrected with a term $\mathbf{c}_k$. The larger the change in parameters $\Delta\boldsymbol{\theta}_k$, the more accurately the model predicts the measured gradients of the optimization problem. On the other hand selecting a large $\Delta\boldsymbol{\theta}_k$ will have a negative effect on the model correction, $\mathbf{c}_k$ (Equation (14)). Since $\mathbf{c}_k$ is based on the linear approximation of the model around $\boldsymbol{\theta}_k$, as $\Delta\boldsymbol{\theta}_k$ increases, the validity of the linear approximation decreases. Thus, as explained in the previous section, to control the accuracy of the linear approximation, a constraint on $\Delta\boldsymbol{\theta}_k$ is imposed by bounding the truncation error $\epsilon_{max}^T$. To summarize, the larger values of $\epsilon_{max}^T$ will allow for large moves in $\Delta\boldsymbol{\theta}_k$ which favours a faster convergence towards the process optimum but this might increase the prediction error incurred by the model. On the other hand, the smaller values will restrict the moves in $\Delta\boldsymbol{\theta}_k$ to generate better predictions but the model may not be able to correct for the optimization gradients accurately, making the algorithm more sensitive to the modelling error. To illustrate this relative effect, the algorithm is solved for $\epsilon_{max}^T = 1\%$ and $5\%$. The convergence in the optimal $S_o$ for these two scenarios is compared in Figure 4. In all the iterations, the input feed rate $F$ converged to the same optimal value (~0.1728 L) so as to satisfy the volume constraint and, therefore, it is not considered further in the discussion. From Figure 4, it is evident that regardless of the choice of $\epsilon_{max}^T$, the algorithm eventually converges to the process optimum. However, the rate of convergence is significantly different for the two cases. In some of the intermediate iterations for $\epsilon_{max}^T = 1\%$, the values of $\Delta\boldsymbol{\theta}_k$ are not sufficiently large for the model to predict the optimization gradients in the correct direction, resulting in an oscillatory and much slower convergence. Whereas, for $\epsilon_{max}^T = 5\%$, the corrections for the optimization are more accurate as a result of which the algorithm converges much faster. However, this improved convergence is at the cost of prediction accuracy. On comparing the total prediction error for all iterations, it was found that on average, this value is nearly 2.5 times higher for $\epsilon_{max}^T = 5\%$.

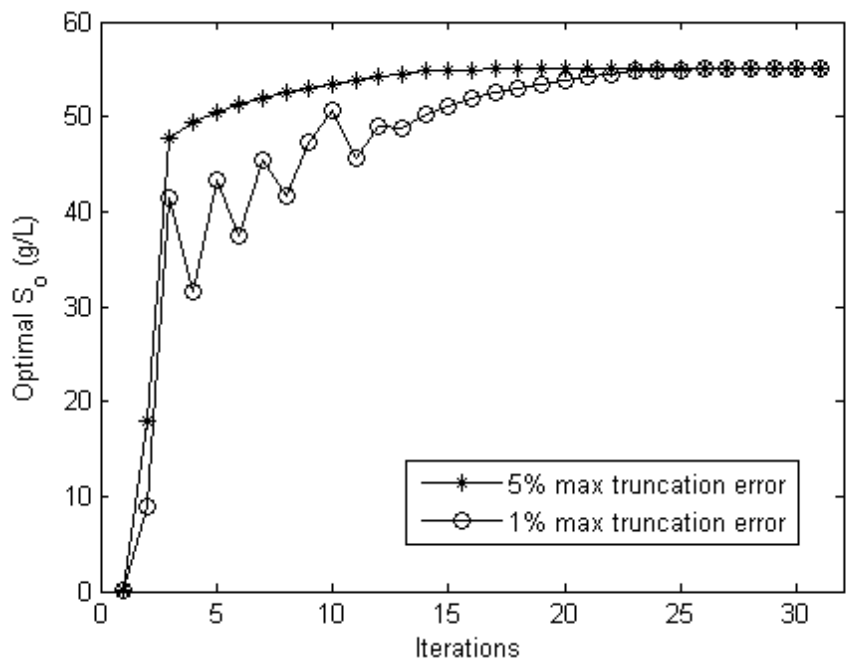

Figure 4: Comparison of the effect of $\epsilon_{max}^T$ on convergence of optimal $S_o$

The key feature of this algorithm, as discussed in previous sections, is not only the convergence to a process optimum but also that the final corrected model predicts the process behaviour accurately and this is corroborated from Figure 5 where the model is used to predict the process variables around the optimum. As can be seen, the predictions are in close agreement with the measurements.

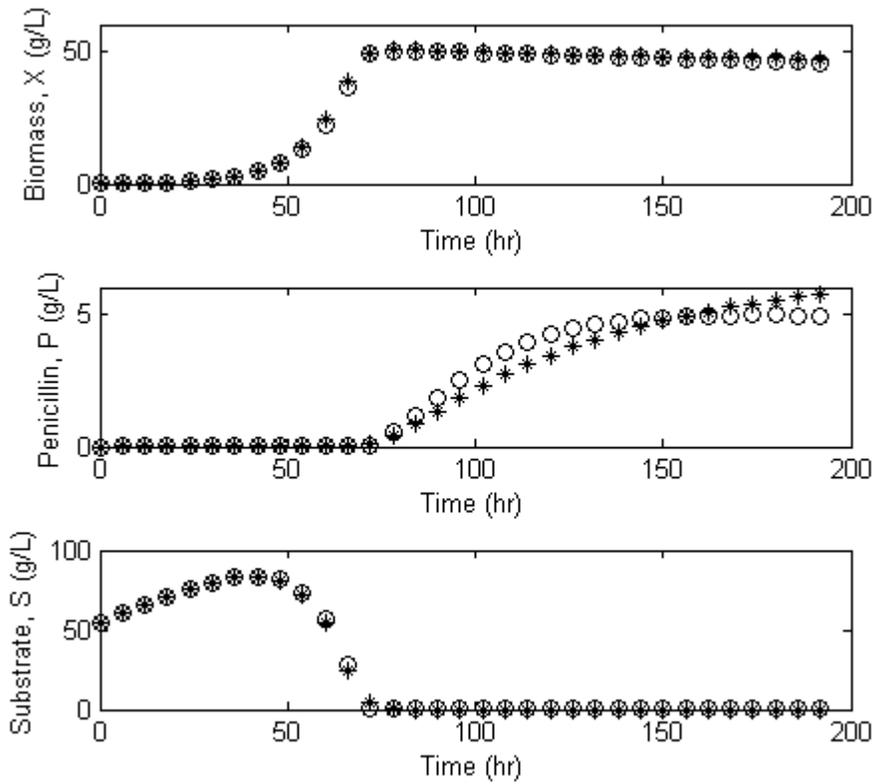

Figure 5: Predictions of corrected model at converged optimal solution

($-*-$ Predictions; $-o-$ Noise-free measurements)

In the next section, the performance of the proposed algorithm is compared with (1) the standard two-step approach and (2) the modifier adaptation algorithms.

*4.2.1. Comparison with Standard "two-step" approach*

The convergence in the optimal $S_o$ corresponding to both "two-step" and proposed methodology is shown in Figure 6. Based on these results, where the proposed algorithm converges to the process optimum $S_o \sim 55$ g/l, the two-step approach converges prematurely to $S_o \sim 34$ g/l. The measured penicillin at the end of batch at $S_o \sim 34$ g/l is $\sim 445$ g which is nearly 25% less as compared to the $\sim 592$g measured at the true process optimum.

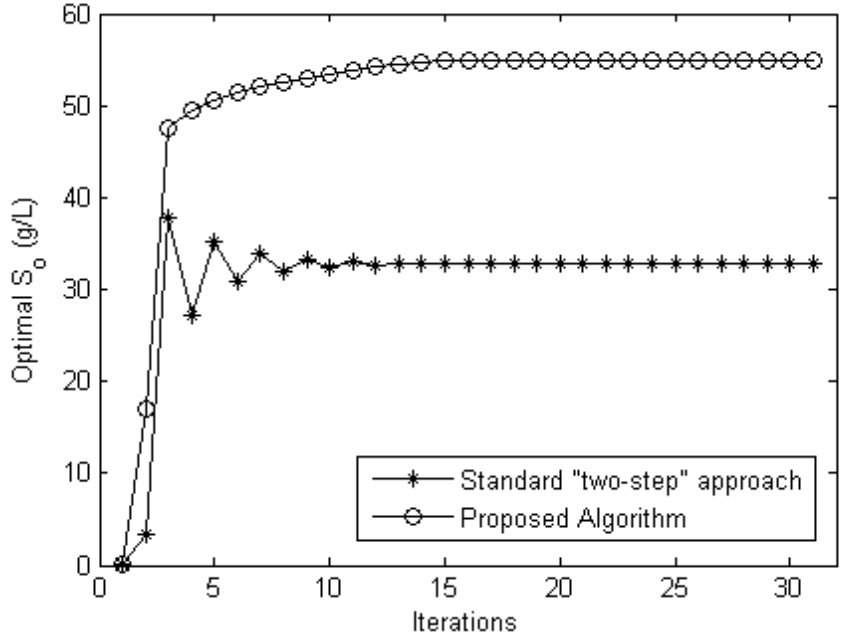

Figure 6: Convergence of optimal $S_o$

The reason for this premature convergence is that the model is inadequate for predicting the true process optimum, or in other words, it cannot satisfy the Equations (2a) and (2b). For all $S_o > 34$ g/l, the change in model parameters, to compensate for the model structure error, is of such an extent that the predicted gradients no longer drive the optimization objective in the direction of its true optimum. On the other hand, in the proposed algorithm, the prediction error and the differences between the predicted and measured gradients of the optimization problem were both used to update the model, thus correcting for the structural uncertainty along the iterations.

*4.2.2. Comparison with modifier adaptation algorithms*

When compared to the class of modifier adaptation algorithms, the proposed correction in this work offers an added advantage in terms of the rate of convergence. As discussed in previous section, in the modifier adaptation algorithms, the optimization objective and the constraints are corrected by adding the differences between the predicted and measured gradients directly to their respective equations (Problem P.4). Since the model is not updated explicitly, the corrected optimization problem may have significant errors in predictions and, generally, this is controlled by filtering the corrections using an empirical filter. For the comparative study, we used first-order exponential filters with three different values for the gain $K = 0.65, 0.5$ and $0.35$. Figure 7, then, compares the convergence in the optimal $S_o$ for the two algorithms for noise free case. In these results, the optimal $S_o$ corresponding to $K = 0.65$ is highly oscillatory and it is quite clear that, without enough filtering, the gradient corrections are more aggressive, leading to significant prediction errors in the input space. However, as the filter gain is decreased, the convergence is much smoother but at the cost of decreased rate of convergence. On the other hand, the proposed algorithm converges much faster and smoothly to the true optimum. The reason is that the model itself is updated to correct for the gradients in the optimization problem which provides a model based filtering that have superior prediction capabilities as compared to the exponential filter in modifier adaptation algorithms.

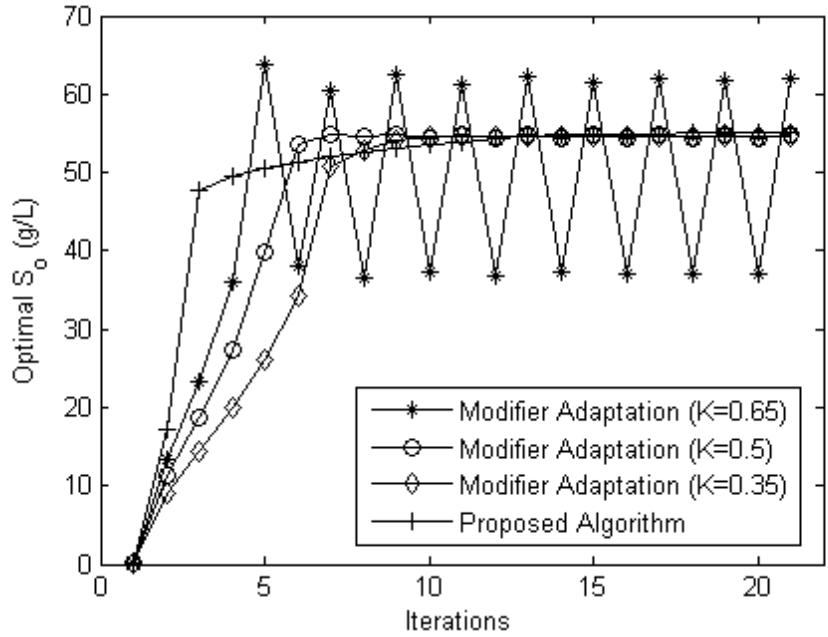

Figure 7: Comparison of proposed and modifier adaption algorithms on convergence of optimal $S_o$

Finally, the convergence of these algorithms is compared in the presence of measurement noise. The modifier adaptation and the proposed algorithms are each solved 10 times with different realizations of the noise and the performance is evaluated in terms of integral absolute error (IAE) and standard deviation in the optimal $S_o$, as summarized in Table 3. The convergence in the average optimal $S_o$ is shown in Figure 8 in the form of error plots. From these results, it is evident that the proposed algorithm is more robust to model errors in the presence of noise. The filter gain that provided smooth convergence in the noise-free situation cannot filter the noise as efficiently as the truncation error ($\epsilon_{max}^T$) in the proposed algorithm. When the corrections are added directly to the optimization problem, the effect of noise in gradients on the optimization objective is additive. Whereas, in the proposed algorithm, the noise in gradients affects the parameter estimates $\boldsymbol{\theta_k} + \Delta\boldsymbol{\theta_k}$ but since the optimization problem is not linear with parameters, this effect is not additive.

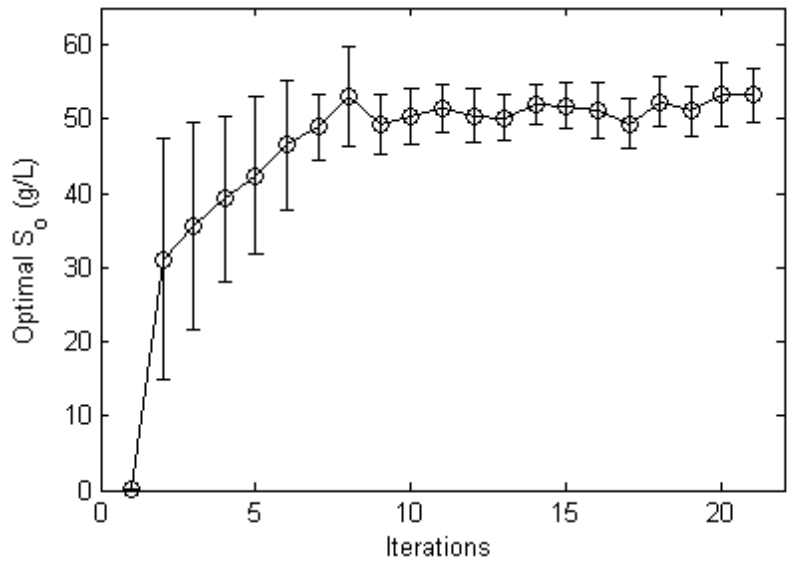

(a)

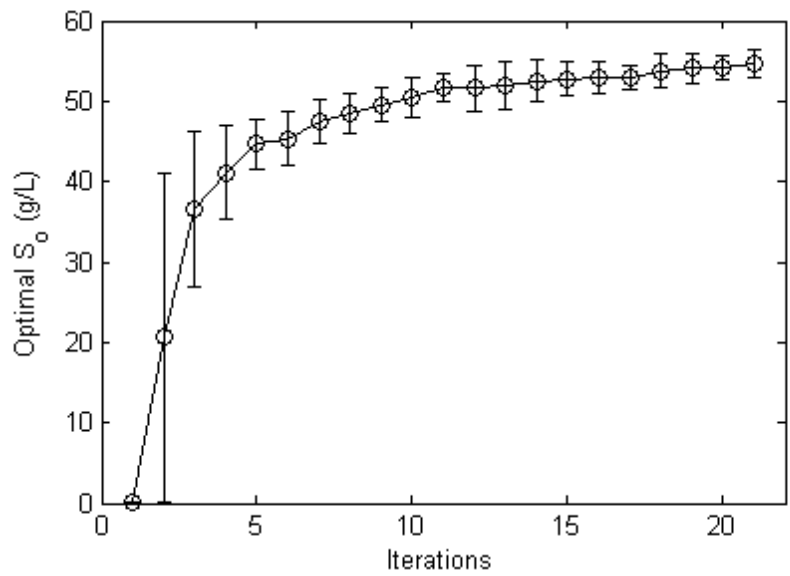

(b)

Figure 8: Average convergence of optimal $S_o$ for (a) Modifier Adaptation algorithm and (b) Proposed algorithm

Table 3: Comparison of convergence properties for the proposed algorithm vs modifier adaptation

|  | IAE | Std. deviation $\sigma$ |
|---|---|---|
| Proposed Algorithm | 8.8597 | 3.5563 |
| Modifier Adaptation | 9.5525 | 5.5805 |

From Figure 8 and Table 3, it can be seen that there is a significant amount of variability in the transient phase for the modifier adaptation algorithm. This is partly related to the fact that the model parameters are never updated, in which case the initial uncertainty in their estimates is propagated throughout the iterations. As a result, for each noise realization, the algorithm may have a significantly different search path if the filter gain is low enough to allow for smaller corrections, as seems to be the case in this example. Increasing the gain decreases this variability in transient but it increases the sensitivity of the algorithm to the noise in gradients, resulting in larger oscillations around the optimum as already shown in the noise free case (Figure 7).

## 5. Conclusions

An iterative optimization algorithm has been proposed where the process models are corrected iteratively for model-plant mismatch in order to guarantee the convergence to the process optimum. The correction is based on linear approximation of the model and is added in a way that upon convergence, the model not only predicts the process behaviour accurately but also satisfies the process optimality conditions. To achieve this goal, the parameter estimation is performed in two sequential steps where a standard estimation problem to minimize the prediction errors is solved first. Then, in the second step, the differences between the predicted and measured gradients of the optimization problem are minimized. The key advantage of this approach is that it provides a model-based filter which has been shown to outperform the exponential filter needed in previously reported algorithms where the gain is selected ad-hoc. The efficiency of the algorithm is illustrated using a fed-batch bioprocess. The rate of convergence depends on the truncation error, used to validate the linear approximation of the correction. For nonlinear models, this approximation may only be valid over a smaller region, therefore limiting the rate of convergence. To this end, an improved approximation of the model has to be used that will be considered in a future study.

## Acknowledgements

The financial support from Natural Science and Engineering Research Council (NSERC), Canada is highly appreciated.

**References**

Bajpai, R.K., Reuss, M. (1980). A Mechanistic Model for Penicillin Production. *Journal of Chemical Technology and Biotechnology*, 30, 332-344.

Beyer, H., Sendhoff, B. (2007). Robust optimization - A comprehensive survey. *Computer Methods in Applied Mechanics and Engineering*, 196 (33-34), 3190-3218

Birol, G., Undey, C., Cinar, A. (2002). A modular simulation package for fed-batch fermentation: penicillin production. *Computers and Chemical Engineering*, 26 (11), 1553-1565.

Chachuat B., Srinivasan, B., Bonvin, D. (2009). Adaptation strategies for real-time optimization. *Computers and Chemical Engineering*, 33 (10), 1557-1567.

Chen, C.Y., Joseph, B. (1987). On-line optimization using a two-phase approach: An application study. *Industrial and Engineering Chemistry Research*, 26, 1924–1930.

Costello, S., Francois, G., Srinivasan, B., Bonvin, D. (2011). Modifier Adaptation for Run-to-Run Optimization of Dynamic Processes. *Proceedings of the 18$^{th}$ IFAC World Congress*, 11471-11476, Milano

Cover T.M, Thomas, J.A. (1991). *Elements of Information theory*, John Wiley and Sons, New York

Diwekar, U.M., Kalagnanam, J.R. (1996). Robust Design using an Efficient Sampling Technique. *Computers and Chemical Engineering*, 20, S389-S394.

Eaton, J.W., Rawlings, J.B. (1990). Feedback control of nonlinear processes using on-line optimization techniques. *Computers and Chemical Engineering*, 14, 469-479

Gao, W., Engell, S. (2005). Iterative set-point optimization of batch chromatography. *Computers and Chemical Engineering*, 29, 1401-1409

Mandur J., Budman, H. (2013). A Robust algorithm for Run-to-run Optimization of Batch Processes, In Proceedings of *10th International Symposium on Dynamics and Control of Process Systems*, 541-546, Mumbai.

Marchetti A., Chachuat, B., Bonvin, D. (2009). Modifier-adaptation methodology for real–time optimization, *Industrial and Engineering Chemistry Research*, 48 (13), 6022-6033.

Marlin, T.E., Hrymak, A.N. (1997). Real-time operations optimization of continuous processes. *In AIChE Symp Ser – CPC – V*, Vol. 93, 156-164.

Martinez, E.C., Cristaldi, M.D., Grau, R.J. (2013). Dynamic optimization of bioreactors using probabilistic tendency models and Bayesian active learning. *Computers and Chemical Engineering*, 49, 37-49

Nagy, Z.K., Braatz, R.D. (2004). Open-loop and closed-loop robust optimal control of batch processes using distributional and worst-case analysis. *Journal of Process Control*, 14(4), 411-422.

Roberts, P.D. (1979). An algorithm for stready-state system optimization and parameter estimation. *International Journal of Systems Science*, 10, 719-734

Ruppen, D., Benthack, C., Bonvin, D. (1995). Optimization of batch reactor operation under parametric uncertainty-computational aspects. *Journal of Process control*, 5(4), 235-240.

Ruppen, D., Bonvin, D., Rippin, D.W.T. (1998). Implementation of adaptive optimal operation for a semi-batch reactor system. *Computers and Chemical Engineering*, 22, 185-199.

Samsatli, N.J., Papageorgiou, L.G., Shah, N. (1998). Robustness metrics for dynamic optimization models under parameter uncertainty. *American Institute of Chemical Engineering Journal*, 44(9), 1993-2006.

Srinivasan, B., Bonvin, D. (2002). Interplay between Identification and Optimization in Run-to-run Optimization Schemes. *In Proc American Control Conference*, 2174-2179, Anchorage

Tatjewski, P (2002). Iterative optimizing set point control-The basic principle redesigned. *In Proceedings of the 15th Triennial IFAC World Congress*

Terwiesch, P., Agarwal, M., Rippin, D.W.T. (1994). Batch unit optimization with imperfect modelling-a survey. *Journal of Process Control*, 4, 238-258